\mag=\magstep1
\input epsf
\documentclass[10pt,a4paper]{amsart}
\usepackage{amssymb,latexsym}

\def\Cal{\mathcal}

\def\mod{\operatorname{mod}}
\def\rank{\operatorname{rank}}
\def\Supp{\operatorname{Supp}}
\def\min{\operatorname{min}}
\def\Res{\operatorname{Res}}

\def\Ind{\operatorname{Ind}}
\def\length{\operatorname{length}}
\def\<<{\langle}
\def\>>{\rangle}

\numberwithin{equation}{section}
\newtheorem{theorem}{Theorem}[section]
\newtheorem{proposition}[theorem]{Proposition}
\newtheorem{corollary}[theorem]{Corollary}
\newtheorem{definition}[theorem]{Definition}

\newtheorem{remark}[theorem]{Remark}
\newtheorem{lemma}[theorem]{Lemma}

\newtheorem{example}[theorem]{Example}

\begin{document}

\title{
Selberg integral and multiple zeta values
}

\author{Tomohide Terasoma}
\address{Department of Mathematical Science, University of Tokyo, Komaba 3-8-1, 
Meguro, Tokyo 153 , Japan}
\email{terasoma@ms.u-tokyo.ac.jp}

%\author{}
%\address{}
%\email{}

%\thanks{The research of the author was supported in part 
%by . 
%}

%    General info
%\subjclass{Primary ; 
% Secondary
% }
\date{\today }
%
%\keywords{}

%%%%%% example %%%%%%
%\begin{equation}
%\label{eq:Omega}
%\end{equation}
%%%%%%%%%%%%%%%%%%%%
%\begin{figure}[t]
%\label{fig:2-move}
%\vskip 10pt
%\centerline{\hbox{\epsfxsize=10cm\epsfbox{s522.eps}}}
%\vskip 5pt
%\caption[]{\label{fig:2-move}
%Graph transformation under a non-trivial 2-move. } 
%\end{figure}
%%%%%%%%%%%%%%%%%%%

\maketitle

\makeatletter
\renewcommand{\@evenhead}{\tiny \thepage \hfill  T.TERASOMA \hfill}
\renewcommand{\@oddhead}{\tiny \hfill  SELBERG INTERAL AND MULTIPLE ZETA FUNCTION
 \hfill \thepage}
\makeatother

\section{Introduction}
\label{sec:introduction}

In this paper, we show that the coefficients of the Talor expansion of
Selberg integrals with respect to its exponent variables are expressed
as a linear combinations of multiple zeta values.  
First object we treat is Selberg integral, the peirod integrals of abelian
coverings of the moduli spaces of $n$-points in $\bold P^1$.
Let $2 \leq r \leq n$ be integers, $\alpha_{ij}$ be positive real numbers.
For an element $f \in \bold C[\frac{1}{x_i -x_j}]$, the function on
$D= \{ x_1 < x_r < \cdots  < x_3 < x_2\}$ defined by
the integral 
$$
\int_{D'} f \prod_{i<j}(x_j -x_i)^{\alpha_{ij}}dx_{r+1} \cdots dx_n,
$$
where $D' =\{ x_1 < x_n < \cdots < x_r\}$ is called a Selberg integral.
It is considered as a family of period integrals for abelian coverings
of the moduli space of distince $n$-points in $\bold C$. It is a function
on $x_1, \dots , x_r$ and exponent parameters $\alpha_{ij}$.
If $r=2$, by the simple functional equation on $x_1$ and $x_2$,
the Selberg integral is determined by its restriction to $x_1 =0$ and
$x_2=1$.
This restricted function
on the exponent parameter $\alpha_{ij}$ 
is called $n-2$-dimensional Selberg integral of $0$-variables.
It is natural to think
that this period integral is equipped with an arithmentic nature.

The second object in our paper is multiple zeta value introduced by Euler.
Let $\bold k =(k_1, \dots ,k_m)$ be a sequence of integers such that
$k_i \geq 1$ $(i=1, \dots m-1)$ and $k_m \geq 2$.  The multiple zeta value
for the index $\bold k$ is define by
$$
\zeta (\bold k) = \sum_{n_1 < \cdots < n_m}\frac{1}{n_1^{k_1}\cdots n_m^{k_m}}.
$$
The natural number $\mid \bold k \mid = \sum_{p=1}^m k_p$ is called the
weight of the index $\bold k$. The integer $\mid \bold k \mid$ is 
called the weight
of the multiple zeta value $\zeta (\bold k)$.  By using the iterated integral 
expression, multiple zeta values are regarded as period integrals for the
fundamental group $\pi_1(\bold P^1 - \{ 0, 1, \infty \})$ of
$\bold P^1 - \{ 0, 1, \infty \}$.  
(See \S \ref{sec:Drinfeld associator} for the iterated integral expression
of multiple zeta values.)
Notice that the motivic weight of 
$\zeta (\bold k)$ is equal to $-2 \mid \bold k \mid$.
The main theorem of this paper is 
\begin{theorem}
For a suitable choice of $f$, the degree $w$ coefficient of the Talor 
expansion for $\alpha_{ij}$ 
of the Selberg integral of $0$-variables
is a linear combination of weight $w$
multiple zeta values.
\end{theorem}
Let us illustrate a primitive example for this statement.  By the well
known equality
$$
\log \Gamma (1 + x) = \gamma x  + \sum_{n \geq 2}\frac{\zeta (n) x^n}{n},
$$
we have
$$
\frac{\Gamma (1+\alpha ) \Gamma (1+\beta )}{\Gamma (1+\alpha + \beta)}
= \exp (\sum_{n \geq 2}\frac{\zeta (n) (\alpha^n+\beta^n - (\alpha+\beta)^n)}{n}).
$$
In this example, we choose $r=2, n=3$ and $f= \alpha /(1-x)$.  
(See \cite{A2} \cite{Z} for another expression of this quantity.)  
We can find the prototype of this theorem in \cite{A2}.
The method
of the choice of $f$ leads us to an interesting combinatorial problem.
In this paper, we will answer to this problem.  This choice of $f$ happens to be
equal to $\beta$-nbc base after Falk-Terao \cite{FT}.

Let us summerize the method of the proof of the main theorem.
Let $\bold C \langle\langle X, Y \rangle\rangle$ be the formal 
non-commtuative free algebra generated by $X$ and $Y$.  After Drinfeld,
associator $\Phi (X,Y)$ is defined as an element of 
$\bold C \langle\langle X, Y \rangle\rangle$. It is known that the
coefficient of $\Phi (X, Y)$ is expressed as a $\bold Q$-linear combination
of multiple zeta values by Le-Murakami \cite{LM}.

Let $\bold Q [[\alpha_i]]$ and $\bold C [[\alpha_i]]$ be  
formal power series rings with vairables
$\alpha_i$ ($i \in I$) over $\bold Q$ and $\bold C$ respectively.
Let 
$\rho : \bold C\langle\langle X, Y \rangle\rangle \to 
M(r, \bold C [[\alpha_i]])$
be a continuous homomorphism, where all the matrix elements of $\rho (X)$ and
$\rho (Y)$ are degree 1 homogeneous in $\bold Q [[\alpha_i]]$.
By the result of Le-Murakami, the coefficients of the Talor expansion 
of any matrix element
$\rho (\Phi (X,Y))$ are expressed by multiple zeta function.
For any solution $s$ of the differential equation
$$
ds = (\frac{\rho (X)}{x}+\frac{\rho (Y)}{x-1})s dx,
$$
we have
$$
\lim_{x \to 1}((1-x)^{-\rho (Y)}s(x)) = \rho (\Phi (X,Y))
\lim_{x \to 0}(x^{-\rho (X)}s(x)).
$$
In this paper, we construct representation $\rho$ and horizontal section $s$
with the following properties.
\begin{enumerate}
\item
All the element of $\lim_{x \to 0}(x^{-\rho (X)}s(x))$
is expressed as $n-3$-dimensional $0$-vairalbe
Selberg integrals by taking a limit for some of $\alpha_i \to 0$.
\item
All the element of $\lim_{x \to 1}((1-x)^{-\rho (Y)}s(x))$
is expressed as $n-2$-dimensional $0$-vairalbe
Selberg integrals by taking the same limit $\alpha_i \to 0$.
\end{enumerate}
For the construction of representation with these properties, we
make a combinatorial preparation in Section \ref{sec:combinatorial preliminaries}.
The construction of $\rho$ depends on the computation of the higher direct 
image of the local system on the moduli space $X_n$ of distinct $n$-points
in $\bold C$ for the projection $X_n \to X_{n-1}$. This computation is
executed in Section \ref{sec:selberg integral}. The limit for $x \to 0$,
$x\to 1$ and $\alpha_i \to 0$ are given in Section 
\ref{sec:proof of the main theorem}.

The author would like to thank T.Kohno for references and discussions.
He also would like to thank M.Kaneko and H.Terao for discussions.

\section{Preliminary}
\label{sec:preliminary}

\subsection{Drinfeld Associator}
\label{sec:Drinfeld associator}
In this section, we  recall known facts about Drinfeld associator. 
Let $R = \bold C \langle\langle  X, Y \rangle\rangle$ 
be the completion of non commutative polynomial
ring of symbols $X, Y$ with respect to its total degree on $X$ and $Y$.
Let $V= R$.  Then $X$ and $Y$ acts on $V$ as the left multiplication
and under this action $X$ and $Y$ are regarded as elements in 
$End_{\bold C}(V)$.  Now we consider a differential form $\omega$
on $\bold C-\{ 0,1\}$ with the coefficient in $End_{\bold C}(V)$ defined
as 
$$
\omega = \frac{X}{x}dx + \frac{Y}{x-1}dx,
$$
where $x$ is the coordinate of $\bold C$.  Let 
$
E(x) = \exp (\int _{x_0}^x \omega)
$
be the solution of the differential equation for $End_{\bold C}(V)$-valued
function $E(x)$
$$
dE(x) = \omega E(x)
$$
with the initial condition $E(x_0) = 1$.
Then by the standard argument for iterated integrals, 
$\exp (\int _{x_0}^x \omega)$ is expressed as
\begin{align}
\label{eq:exponential}
\exp (\int _{x_0}^x \omega) = 1 + \int _{x_0}^x \omega +
\int _{x_0}^x \omega\omega + \cdots .
\end{align}
Here we use the convention for iterated integrals defined by
the inductive relation
$$
\int_p^q \omega_1 \cdots \omega_n =
\int_p^q (\omega_1(q_1)\int_p^{q_1}\omega_2 \cdots \omega_n).
$$
The expression ($\ref{eq:exponential}$) implies, $\exp (\int _{x_0}^x \omega) 
\in \bold C \langle\langle X, Y \rangle\rangle ^{\times}$ and the shuffle relation for
iterated integral imples that
$E = \exp (\int _{x_0}^x \omega)$ is a group like element, i.e.
$\Delta (E) = E \otimes E $ in $\bold C\langle\langle X, Y\rangle\rangle  \hat\otimes
\bold C \langle\langle X, Y\rangle\rangle $ where the comultipication $\Delta :\bold C\langle\langle X, Y\rangle\rangle  \to
\bold C\langle\langle X, Y\rangle\rangle  \hat\otimes \bold C \langle\langle X, Y\rangle\rangle $ is given by $\Delta (X)
= X\otimes  1 + 1 \otimes X$ and $\Delta (Y) =Y \otimes 1 + 1 \otimes Y$.
The set $\hat G = \{ \Delta (g) = \Delta (g) \otimes \Delta (g) \mid g 
\in \bold C\langle\langle X, Y\rangle\rangle ^{\times}\}$ is called the set of group like elements
and closed under the multiplication.  By the theory of differential equation
only with regular singularity, the limit
$
\lim_{x\to 1} \exp (\int_x^0 \frac{Y}{x-1}dx)
\exp (\int _{x_0}^x \omega) 
$
exists. In the same way, the limit
$$
\Phi (X, Y) = \lim_{x\to 1, y\to 0} \exp (\int_x^0 \frac{Y}{x-1}dx)
\exp (\int _{y}^x \omega) \exp (\int_1^y \frac{X}{x}dx)
$$
exists and contained in $\bold C\langle\langle X, Y\rangle\rangle ^{\times}$. 
$\Phi (X, Y)$ is called the Drinfeld associator. Since 
$\exp (\int_x^0 \frac{Y}{x-1}dx)$ and $\exp (\int_1^y \frac{X}{x}dx)$
are elements in $\hat G$ and $\hat G$ is a closed subset of 
$\bold C\langle\langle X, Y\rangle\rangle ^{\times}$,
the limit $\Phi (X,Y)$ is an element in $\hat G$.

We recall that the relation between multiple zeta values and the
coefficients of the Drinfeld associator.  Firstly, we recall the definition
of multiple zeta values. Let $k_1, \dots , k_n$ be integers such that
$k_i \geq 1$ for $i = 1, \dots ,n$ and $k_n \geq 2$. Set 
$\bold k = (k_1, \dots , k_n)$.  The following series 
$$
\zeta (\bold k) = \zeta (k_1, \dots, k_n)
=\sum_{m_1 < m_2 < \cdots < m_n}\frac{1}{m_1^{k_1}m_2^{k_2}\cdots m_n^{k_n}}
$$
is called the multiple zeta values for the index $\bold k = (k_1, \dots ,k_n)$.
The number $\mid \bold k \mid = \sum_{i=1}^n k_i$ is called the weight of 
the index $\bold k$. Let $L_w$ be the finite dimensional $\bold Q$
vector subspace of $\bold C$ generated by $\zeta (\bold k)$, with
the weight $\mid \bold k \mid = w$.  The following iterated
integral expression of multiple zeta values is fundamental.
\begin{align*}
\zeta (k_1, \dots , k_n)  = & \int_0^1
\underbrace{\frac{dx}{x}\cdots \frac{dx}{x}}_{k_n -1}\frac{dx}{1-x}
\underbrace{\frac{dx}{x}\cdots \frac{dx}{x}}_{k_{n-1} -1}\frac{dx}{1-x} \\
& \cdots
\underbrace{\frac{dx}{x}\cdots \frac{dx}{x}}_{k_1 -1}\frac{dx}{1-x}.
\end{align*}
By using this expression and shuffle relation, 
for elements $a$ and $b$ in $L_{w_1}$ and $L_{w_2}$, 
we can show that $ab$ is an element
in $L_{w_1 + w_2}$.
Using this fact, we define the homogeneous multiple zeta value ring 
(homogeneous MZV ring for short) $H$ in
$\bold C\langle\langle X, Y\rangle\rangle $ by
$$
H= \oplus_{w \geq 0}\oplus_{W:\text{word of length $w$ on $X,Y$}}
L_w\cdot W.
$$
The following proposition is due to Le-Murakami \cite{LM}.
\begin{proposition}
\label{prop:L-M}
$\Phi (X, Y) \in H$.
\end{proposition}
It is very usuful to specialize this universal result to a special class of
representations of $\bold C \langle\langle X, Y \rangle\rangle$.
Let $R$ be a homogeneous complete ring generated by degree 1 elements
over $\bold Q$,
i.e. R is generated topologically by degree 1
homogeneous elements $\alpha_1, \dots ,\alpha_m$
with homogeneous relations and complete under the topology defined by its
degree.  The decomposition of R with respect to its degree is denoted by
$R = \hat\oplus_{d \geq 0} R_d$.  Let $\bold R_{\bold C}$ be the completion of
$R\otimes \bold C$ with resepct to the topology defined by its degree.
A ring homomorphism $\rho :\bold C \langle\langle X, Y\rangle\rangle  \to M(r, R_{\bold C})$
is called a homogeneous rational representation of degree 1 if and only
if all the matrix elements of $\rho (X)$ and $\rho (Y)$ are degree
1 homogeneous elements in $R$.  The homogeneous MZV 
ring $H_R$ for $R$ is defined by 
$H_R = \hat\oplus_{d \geq 0}(R_d \otimes L_d)$.
The following corollary is a direct consequence of Proposition \ref{prop:L-M}.
\begin{corollary}
\label{corollary:homogeneous rational representation}
Let $\rho : \bold C \langle\langle X, Y\rangle\rangle  \to M(r, R_{\bold C})$ be a homogeneous
rational representation of degree 1.  Then all the matrix elements 
of $\rho (\Phi (X,Y))$
are elements in $H_R$.
\end{corollary}

\section{Selberg integral}
\label{sec:selberg integral}
\subsection{Combinatorial aspects}
\label{subsec:combinatorial aspects}
Let $[n] =\{ 1,\dots ,n\}$.  A graph $\Gamma$ consisits of the set of vertices 
$V_{\Gamma}$ and edges $E_{\Gamma}$. We assume
that every edges have distinct two terminals.
Moreover we assume for any two vertices $p$ and $q$, there exists at most
one edge whose terminals are $p$ and $q$. An edge is written
as $(p,q)$, where $p$ and $q$ are its terminals.
For a graph $\Gamma$, we can associate a 
1-dimensional simplicial complex 
by usual manner and we use the standard terminology connected 
component, tree, and so on. Moreover, if the
order of $E_{\Gamma}$ is specified, it is called ordered graph. A specified point
for each connected components is called a root.  
The specified point for a connected component is called the root of 
the connected component and the set of roots is denoted by
$R= R_{\Gamma}$.  
For two sets $V$, $R$ such that $R \subset V$, we define
$\Omega^i (V\mod R)$ by $\wedge^i(\Omega_{X_V}^1/p^*\Omega_{X_R}^1)$,
where $X_V = \{ (x_i)_{i \in V} \mid x_i \neq x_j \text{ for }
i \neq j\}$,
$X_R = \{ (x_i)_{i \in R} \mid x_i \neq x_j \text{ for }
i \neq j\}$, and $p$ is the natural projection $X_V \to X_R$. 
Then it is easy to see that $\Omega^{\# V -\# R} (V\mod R)$ is a rank 1 
$\Cal O_{X_n}$ module generated by $\wedge_{i \in V - R} dx_i$.
For an edge $e = (p,q) , p,q \in V_{\Gamma}$, we define 
$\omega_e = d\log (x_p - x_q) \in \Omega^1(V \mod R)$.  
For an ordered tree, we define 
$\omega_{\Gamma}$ as
$\omega_{\Gamma} = \omega_{e_r}\wedge
\cdots \wedge \omega_{e_1}$ in $\Omega(V \mod R)$, where 
$E_{\Gamma}=\{ e_1, \dots ,e_r \}$ and $e_1 < \cdots < e_r$.
It is easy to to see the following lemma.

\begin{lemma}
Assume $\# E = \# V -\#R$. Then
$\Gamma$ is a tree if and only if $\omega_{\Gamma} \neq 0$
\end{lemma}

Let $R$ be a sub set of $[n]$ such that 
$\{ 1, 2\}\subset R$.
We define an ordering of $[n]$ by
$1 \ll n \ll \cdots \ll 3 \ll 2$.
We define $D(R)$ by  
$\{ (x_1, \dots ,x_i)_{i \in R}\mid x_i <x_j \text{ for }
i\ll j \}$. For two subsets $V$ and $R$ of $[n]$ such that $R \subset V$, 
the fiber of the map $D(V) \to D(R)$ at
$(x_i)_{i \in R} \in D(R)$ is denoted by 
$D(V/R,x_i)_{i \in R}$.
Let $\alpha_{i,j}$ $(i,j \in V)$ be positive real numbers. 
We choose a branch of
$\Phi (V) = \prod_{i \ll j}(x_j - x_i)^{\alpha_{i,j}}$ 
on $D(V)$ with $\Phi \in \bold R_+$.
For an ordered rooted graph $\Gamma$ whose root set is $R$, we define 
a funciont $S_{\Gamma }=S_{\Gamma}(V/R, x_i)_{i \in R}$ on $D(R)$ as
$$
S_{\Gamma }(V/R, x_i)_{i \in R} = \int_{D(V/R,x_i)_{i \in R}}
\Phi (V)\prod_{(i,j) \in E_{\Gamma}} 
\alpha_{i,j}\omega_{\Gamma}
$$
If $R$ is fixed, it is denoted by $S_{\Gamma}$.  Then $S_{\Gamma}$
is a function on $(x_i)_{i \in R}$ and $\alpha_{i,j}$.  The free abelian
group generated by ordered rooted graphs whose root set and the set of vertices
are $R$ and $V$, is denoted by $\Gamma (V,R)$.  For an element
$\gamma = \sum a_{\Gamma} \Gamma$ in $\Gamma (V,R)$, we define $S_{\gamma}$
by $S_{\gamma}= \sum a_{\Gamma}S_{\Gamma}$. The function
$S_{\gamma}$ is called the Selberg integral for $\gamma$.

Before the presentation of the main theorem, 
we introduce several combinatorial notions.
For two natrual numbers $n,r$ such that $2 \leq r \leq n$, we set
$R= [r]$ and $V = [n]$. For an ordered rooted graph $\Gamma$,
whose vertex set and root set are $[n]$ and $[r]$, we define an element
$\Gamma \wedge (n+1, i)$ in $\Gamma ([n+1],[r])$ for $i \in [n]$
according to the following recipe.
\begin{enumerate}
\item
Choose a subset $A$ of edges which are connecting to $i$.
($A$ may be an emptyset.)
\item
Replace the number $i$ by $n+1$ for all the edges contained in $A$ choosen in
1.
\item
Make a graph $\Gamma_A$ by adding the edge $(n+1, i)$ to the graph $\Gamma$
and extend the original ordering to that of the edge set of $\Gamma_A$
such that $(n+1,i)$ is the biggest edge.
\item
Consider the sum $\sum_A \Gamma_A$ of $\Gamma_A$, where $A$ runs through
all the subsets of edges connecting to $i$.  
This summation
is denoted by $\Gamma \wedge(n+1, i)$.
\end{enumerate}
We extend the operation $\wedge (n+1, i_{n+1})$ from $\Gamma ([n], [r])$
to $\Gamma ([n+1], [r])$ by linearity.
For an element $\gamma \in \Gamma ([l], [r])$ and 
$(l+1, i_{l+1}),\dots ,(n,i_n)$, where $i_{l+1} \in [l], \dots ,i_n \in [n-1]$,
we define 
$\gamma\wedge (l+1, i_{l+1})\wedge\cdots \wedge (n, i_n)$ inductively by
$$
\gamma\wedge (l+1, i_{l+1})\wedge \cdots \wedge (n, i_n) =
(\gamma\wedge (l+1, i_{l+1})\wedge \cdots 
\wedge (n-1, i_{n-1}))\wedge (n, i_{n})
$$
The graph $\Gamma$ with $V_{\Gamma} = R$ and $E_{\Gamma} = \emptyset$
is denoted by $\emptyset (R)$. A grpah is denoted by $e_1 e_2 \cdots  e_b$,
where the set of edges is 
$\{e_1 < e_2 < \cdots  < e_b \}$ .
\begin{example}
If $R =\{1,2\}, i_3 = 2, i_4=2$, then
$$
\emptyset (R)\wedge (3,2)\wedge (4,2) = (2,3)\wedge (4,2) = (2,3)(4,2) + (4,3)(4,2)
$$
\end{example} 
We state the main theorem.  Let $H_{\alpha}$ be a homogeneous MZV ring for 
$\bold Q \langle\langle \alpha_{i,j},\alpha_{1,k}, 
\alpha_{2,k}\rangle\rangle _{3 \leq i,j,k \leq n,
i\neq j}$.
\begin{theorem}
\label{theorem:the main theorem}
Let $R=\{1,2\}$. For any $i_3 \in [2], \dots ,
i_n \in [n-1]$, put $\gamma = \emptyset (R)\wedge (3, i_3) \wedge \cdots \wedge
(n, i_n)$. Then
$S_{\gamma}([n]/[2],0, 1)$ is a holomorphic 
function on $\alpha_{i,j}$ and an element of $H_\alpha$. 
\end{theorem}

\subsection{Differential equation satisfied by Selberg integral}
\label{subsec:differential equation}

First we compute the higher direct image of the connection on the
configuration space $X_n =\{(x_1, \dots ,x_n) \mid x_i \neq x_j 
\text{ for $i \neq j$ } \}$
for distinct $n$-points in $\bold C$ for the morphism
$\pi: X_n \to X_{n-1}$ defined by $(x_1, \dots , x_n) \to (x_1, \dots ,x_{n-1})$.
Let $A_{i,j} \in M(d,\bold C)$ be matrices for $1 \leq i \neq j \leq n$
satisfying the following relations.  These realations are called
pure braid relations.
\begin{enumerate}
\item
$A_{i,j} = A_{j,i}$.
\item
$[A_{i,j}, A_{k,l}]=0$ for all distinct $i, j, k, l$.
\item
$[A_{i,j} + A_{j,k}, A_{i,k}] = 0$ for all distinct $i, j, k$.
\end{enumerate}
Then the matrix valued 1-form
$$
\omega = \sum_{1 \leq i<j \leq n} A_{i,j}d \log (x_i - x_j)
$$
defines an integrable connection $\nabla$ on 
$\Cal O_{X_n}^d = \{ v= ^t(v_1, \dots ,v_d)\}$
by
$$
\nabla v = dv - \omega v.
$$
Let $v$ be a horizontal section of
the connection $\nabla$ on $D([n])$, i.e. $dv = \omega v$.  For 
$i \in [n-1]$, and $(x_1, \dots ,x_{n-1})\in D([n-1])$, 
we define $w_i$ as
$$
w_i = \int_{D([n]/[n-1],x_1, \dots ,x_{n-1})}\frac{A_{n,i}}{x_n - x_{i}}v dx_n.
$$
Then $w_i$ is a function on $(x_1, \dots , x_{n-1})\in D([n-1])$.
We have the following proposition.
\begin{proposition}
\label{proposition:higher direct image}
\begin{enumerate}
\item
$w_1 + \cdots + w_{n-1} = 0$.
\item
Let $W = ^t(w_1, \dots , w_{n-1})$. Then $W$ satisfies the differential equation
$$
dW = \sum_{1\leq i<j \leq n-1}\frac{A'_{i,j}(dx_i - dx_j)}{x_i - x_j}W,
$$
where
\begin{align}
\label{eq:big matrix}
A'_{i,j} =
\left(\begin{matrix}
A_{ij} & \overset{i}{\cdots} & \overset{j}{\cdots} & 0 \\
\vdots & A_{ij}+A_{nj} & -A_{ni} & \vdots \\
\vdots & -A_{nj} & A_{ij}+A_{ni} & \vdots \\
0 & {\cdots} & {\cdots} &A_{ij}  \\
\end{matrix}\right)
\begin{matrix}
 \\ i \\ j \\ \\
\end{matrix}
\end{align}
\end{enumerate}
\end{proposition}
\begin{proof}
1.
By the equality 
$$
\frac{\partial v}{\partial x_n} = \sum_{j=1}^{n-1}\frac{A_{nj}}{x_n-x_j} v,
$$
and Stokes' theorem, the equatity follows.

2.
By using the differential equation for $v$, we have
\begin{align*}
& \frac{\partial}{\partial x_i}(\frac{A_{nj}}{x_n- x_j}v) \\
= &
\frac{A_{nj}}{x_n- x_j}(\sum_{k \neq i,j,n}\frac{A_{ik}v}{x_i- x_k}
+\frac{A_{ij}v}{x_i- x_j}+\frac{A_{ni}v}{x_i- x_n}) \\
= &
\sum_{k \neq i,j,n}\frac{A_{ik}A_{nj}v}{(x_i- x_k)(x_n- x_j)}+
\frac{A_{nj}v}{x_i- x_j}\{
-\frac{A_{ni}v}{x_n- x_i}+\frac{(A_{ij}+A_{ni})v}{x_n- x_j} 
\} .
\end{align*}
By the commutativity condition of $A_{ij}$, we have
$$
\frac{\partial}{\partial x_i} w_j = \sum_{k \neq i,j, 1 \leq k \leq n-1}
\frac{A_{ik}}{x_i - x_k}w_k + \frac{1}{x_i- x_j}
\{ (A_{ij}+ A_{ni}) w_j - A_{nj} w_i\}
$$
for $i \neq j$. Using the relation in 1., $\frac{\partial}{\partial x_i}w_i =
\sum_{j \neq i, 1 \leq j \leq n-1}\frac{\partial}{\partial x_i}w_j$,
and we obtain the statement of 2.
\end{proof}
\begin{remark}
\label{remark:homomorphism}
If $A_{ij}$ satisfies the infinitesimal pure braid relations, then
matrices $A_{ij}'$ defined by (\ref{eq:big matrix}) 
also satisfies the inifinitesimal
pure braid relations for $n-1$.  Therefore the connection $\nabla '$
given by 
$$
\nabla ' W = dw -\sum_{1 \leq i <j \leq n-1}\frac{(dx_i -dx_j)A_{ij}'}{x_i -x_j} W
$$
on $(\Cal O^{\oplus d})^{\oplus (n-1)}$ is integralble.
\end{remark}

Let $V_n$ be the local system of horizontal sections of the connection $\nabla$
and $V_n\mid_{ \pi^{-1}(x_1^0, \dots , x_{n-1}^0)}$ its restriction of the fiber
$\pi^{-1}(x_1^0, \dots , x_{n-1}^0)$ of 
$\pi :X_n \to X_{n-1}$.  Then the Eular-Poincare
characteristic of $V_n\mid_{\pi^{-1}(x_1^0, \dots , x_{n-1}^0)}$ is
$-(\rank V)\cdot (n-2)$.  Therefore under certain non-resonance condition,
$\dim H^1(\pi^{-1}(x_1^0, \dots ,x_{n-1}^0), V_n)$ is equal to
$\rank V\cdot (n-2)$.  By a direct comupation, the submodule 
$(\Cal M )^{red} = 
\{ W = ^t(w_1, \dots ,w_{n-1}) \mid w_i \in \Cal O^{\oplus d},
\sum_{i=1}^{n-1} w_i = 0\}$ comes to be a sub connection of 
$\Cal M =( (\Cal O^{\oplus d})^{\oplus (n-1)}, \nabla ')$.
As a consequence, horizontal section of $(\Cal M )^{red}$ is equal to the
higher direct image of $V$ under the projection $\pi$.  This construction is
compatible with the sub local system in $V$.

  We apply this inductive formula to compute the differential equation satisfied
by Selberg integral.  Note that the 
similar computation is executed in \cite{A1}
with a different base of de Rham cohomology.  
Selberg integral is holomorphic with respect to  $\alpha_{ij}$ for our
base. 
This base is nothing but the $\beta$-nbc base introduced in
Falk-Terao \cite{FT}.

Let $R$ be a ring. For a set of elements 
$\bold a =\{ a_{pq}\}_{1 \leq p< q \leq k}$ satisfying the infinitesimal pure
braid relation, we define a set of elements 
$\Ind( \bold  a) = \{ \Ind (\bold a)_{ij} \}_{1 \leq i < j \leq k-1}$
in $M(k-1,R)$ by
$$
\Ind(\bold a)_{ij}= \left(
\begin{matrix} 
a_{ij} & \overset{i}{\cdots} & \overset{j}{\cdots} & 0 \\
\vdots & a_{ij}+a_{kj} & -a_{ki} & \vdots \\
\vdots & -a_{kj} & a_{ij}+a_{ki} & \vdots \\
0 & {\cdots} & {\cdots} &a_{ij}  \\
\end{matrix}\right)
\begin{matrix}
\\ i \\ j \\ \\
\end{matrix}.
$$
Let $2 \leq r \leq n$ be integers and $V_{r,n}$ be a $\bold C$ vector
space of dimension $r(r+1)\cdots (n-1)$ whose coordinates are given by
$v_{i_{r+1}, \dots ,i_n}$ for 
$1 \leq i_{r+1} \leq r, \dots , 1\leq i_n \leq n-1$.
We define ${\bold A}^{(p)} =\{ A^{(p)}_{ij}\}_{1\leq i < j \leq p} $ 
for $p= r, \dots , n-1$ by
$$
{\bold A}^{(p)}=\Ind ({\bold A}^{(p+1)}) 
$$
and $A^{(n)}_{ij} = a_{ij}$.  
We define $V_{k,n}$
valued function $S^{(k)}(x_1, \dots , x_k)$ on $D([k])$ inductively by
$$
S^{(k)}=
\left(\begin{matrix}
\int_{D([k+1]/[k],x_i)_{i \in [k]}}\frac{A^{(k+1)}_{k+1,1}}{x_{k+1}- x_1}
S^{(k+1)}(x_1, \dots ,x_{k+1})dx_{k+1} \\
\vdots \\
\int_{D([k+1]/[k],x_i)_{i \in [k]}}\frac{A^{(k+1)}_{k+1,k}}{x_{k+1}- x_k}
S^{(k+1)}(x_1, \dots ,x_{k+1})dx_{k+1} \\
\end{matrix}\right)
$$
for $k = r,\dots ,n-1$ and 
$$
S^{(n)} = \prod_{1\leq i \ll j \leq n}(x_j - x_i)^{\alpha_{ij}}.
$$
We have the following corollary of Proposition 
\ref{proposition:higher direct image}.
\begin{corollary}
The $V_{k,n}$ valued function $S^{(k)}$ 
satisfies the following differential equation
$$
dS^{(k)} = \Omega_k S^{(k)},
$$
where $\Omega_k = \sum_{1 \leq i < j \leq k}\frac{A^{(k)}_{ij}d(x_i -x_j)}
{x_i - x_j}$.
\end{corollary}
The next proposition is used at the proof of Main theorem 
\ref{theorem:the main theorem}.
\begin{proposition}
\label{proposition:sum relation}
Let $S^{(r)}_{i_{r+1}, \dots ,i_n}$ be the $(i_{r+1}, \dots ,i_n)$-component 
of $S^{(r)}$.
Then we have
\begin{align}
\label{eq:fundamental relation}
\sum_{i_p' =1}^{p-1}
S_{i_{r+1}, \dots , ,i_{p-1}, i_p', i_{p+1}, \dots , i_n} = 0.
\end{align}
\end{proposition}
\begin{proof}
If $p= r+1$, then it is nothing but the first statement of 
Proposition \ref{proposition:higher direct image}. Suppose $p > r+1$.
Then $(i_{r+1}, \dots ,i_{p-1})$-part of $S^{(r)}$ is a linear combination
of 
\begin{align}
\label{eq:sum relation}
 \int \prod_{i=r+1}^{p-1}A_{p_iq_i}^{(p)}
\prod_{j=r+1}^{p-1}\frac{1}{x_j- x_{i_j}}
S^{(p)}dx_{r+1}\cdots dx_{p-1} 
\end{align}
Since the set $\{(a_{i_p, \dots ,i_n}) \mid 
\sum_{i_p'=1}^{p-1} a_{i_p', \dots ,i_n} = 0 \}$ is stable under the action of
$A_{ab}^{(p)}$, (\ref{eq:sum relation}) satisfies the relation
$\sum_{i_p'=1}^{p-1} a_{i_p', \dots ,i_n} = 0$.
\end{proof}

\section{Combinatorial Preliminaries}
\label{sec:combinatorial preliminaries}

\subsection{Statement of the main theorem}
\label{subsec:statement of the main theorem}

In this section, we present combinatorial facts which is used to the
computataion of Selberg integrals.  Let $P_n$ be the non-commutative ring
$\bold C [a_{ij}]$ with the generators $a_{ij}$ $(1 \leq i,j \leq n)$
and the infinitesimal pure braind relations.
We define a set of matrices 
$\bold A^{(k)} = \{ A_{ij}^{(k)}\}_{1 \leq i , j \leq  k}$
in $M(k(k+1) \cdots (n-1), P_n)$ inductively by the relations:
$$
{\bold A}^{(k)} = \Ind (\bold A^{(k+1)})
$$
for $k = r, \dots , n-1$ and $A_{ij}^{(n)} = a_{ij}$.  We introduce the degree
of $P_n$ by $\deg a_{ij} = 1$.  Then the matrix elements of $A_{ij}^{(k)}$ are
degree $1$ for $k = r, \dots ,n$ and $A_{ij}^{(k)}$ satisfes the pure braid
relations.
In other words, 
A ring homomorphism 
$P_r \to M(r(r+1)\cdots (n-1), P_n)$ is defined by attaching
$A_{ij}^{(r)}$ to $a_{ij}\in P_r$.  We define a vector 
$w_k \in P_n^{k(k+1)\cdots (n-1)}\otimes \bold C[\frac{1}{x_i -x_j}]$
inductively by the relation
\begin{align}
\label{eq:tensor vector}
w_k=
\left(\begin{matrix}
\frac{A^{(k+1)}_{k+1,1}}{x_{k+1}- x_1}
w_{k+1} \\
\vdots \\
\frac{A^{(k+1)}_{k+1,k}}{x_{k+1}- x_k}
w_{k+1} \\
\end{matrix}\right)
\end{align}
for $k = r, \dots, n-2$ and
$$
w_{n-1}=
\left(\begin{matrix}
\frac{A^{(n)}_{n,1}}{x_{n}- x_1} \\
\vdots \\
\frac{A^{(n)}_{n,n-1}}{x_{n}- x_{n-1}} \\
\end{matrix}\right).
$$
In this section, we express each coordinate of $w_r$ in terms of
combinatorics introduced in \S \ref{subsec:combinatorial aspects}.
For an ordered tree $\Gamma$ with the vertex set 
$[k]$ and root set $R = [r]$,
we define $A_{\Gamma}^{(k)}$ by
$$
A_{\Gamma}^{(k)} = \prod_{i=l}^1 A_{p_i,q_i}^{(k)} \in M(k(k+1)\cdots (n-1), P_n),
$$
where $E_{\Gamma} =\{ e_1 < \cdots < e_l\}$ and $e_i = (p_i, q_i)$.
Here we use the notation
$\prod_{i=l}^1 a_i = a_l a_{l-1}\cdots a_1$ in a non-comutative
ring.  We define a matrix valued differential form
$\eta_{\Gamma} \in \Omega ([k]\mod [r]) \otimes M(k(k+1) \cdots (n-1), P_n)$
by
$\eta_{\Gamma} = A_{\Gamma}^{(k)}\omega_{\Gamma}$,
where $\omega_{\Gamma}$ is defined in \S \ref{subsec:combinatorial aspects}.
For an element $\gamma \in \Gamma ([k], [r])$,
we define $\eta_{\gamma}$ by
$\eta_{\gamma} = \sum_{\Gamma} a_{\Gamma} \eta_{\Gamma}$,
where $\gamma = \sum_{\Gamma} a_{\Gamma}\Gamma$.
\begin{theorem}
\label{theorem:combinatorics}
Let us denote the $(i_{r+1}, \dots ,i_n)$-th coordinate of $w_r$
by 
\linebreak
$w_r(i_{r+1}, \dots , i_n)$.  Then 
$$
w_r(i_{r+1}, \dots , i_n)dx_n\wedge \cdots \wedge dx_{r+1} =
\eta_{\gamma},
$$
where $\gamma = \emptyset ([r]) \wedge (r+1, i_{r+1}) \wedge \cdots
\wedge (n, i_n)$.
\end{theorem}
The rest of this section is spent to prove 
Theorem \ref{theorem:combinatorics}.

\subsection{Several lemmata}
\label{several lemmata}

Let $\Gamma$ be an ordered graph with  the root set $[r]$ and the vertex
set $[n-1]$.  The edge set is denoted by $E=\{ e_1 < \cdots < e_l\}$
and $e_i = (p_i, q_i)$.
Suppose that $p$ and $q$ are contained in the
same connected component.  Then there exists unique path $P$ connecting $p$
and $q$ in $\Gamma$.  
We write $P = \{ e_{t_1}, \dots , e_{t_m}\}$. The subgraph
$P$ looks like figure 1.
\begin{figure}
%\label{fig:path}
\vskip 10pt
\centerline{\epsfbox{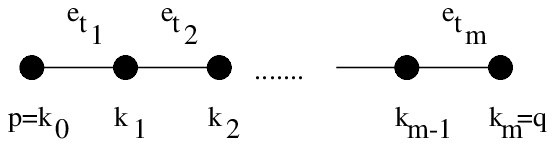}}
\vskip 5pt
\caption[]{ } 
\end{figure}
\begin{lemma}
\label{lemma:product}
Let $A_{\Gamma}^{(n-1)} \in M(n-1, P_n)$ be defined as in \S 
\ref{subsec:statement of the main theorem} 
\begin{enumerate}
\item
If $q$-th component of 
\begin{align}
A_{\Gamma}^{(n-1)}\left(
\begin{matrix}
0 \\ \vdots \\ a_{np} \\ \vdots \\ 0 
\end{matrix} \right) \in P_n^{(n-1)}
\end{align}
is not zero, then $p$ and $q$ are contained in the same connected
component and $t_1 < t_2 < \cdots < t_m$.
\item
Suppose that $t_1 < t_2 < \cdots < t_m$.  We write vertices of the
path $P$ as $p= k_0, k_1, \dots , k_m = q$, (see figure 1)
and define $B_i$ $(i=1, \dots ,l)$ by
$$
B_i = \begin{cases}
-a_{k_j,n} & (\text{ if } i = t_j) \\
a_{p_iq_i} + a_{n q_i} & (\text{ if }  t_j < i < t_{j+1} \text{ and  $e_i$
ajacents to $k_j$ and put $p_i=k_j$}) \\
a_{p_iq_i}  & (\text{ if }  t_j < i < t_{j+1} \text{ and  $e_i$
does not ajacent to $k_j$}) \\
\end{cases}
$$
(For the second case see figure 2.)
\begin{figure}
%\label{fig:path}
\vskip 10pt
\centerline{\epsfbox{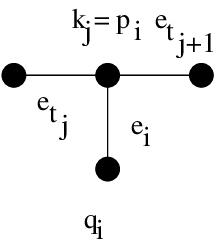}}
\vskip 5pt
\caption[]{ } 
\end{figure}
Then $q$-th component of (1) is equal to $\prod_{i=l}^1 B_i$.
\end{enumerate}
\end{lemma}
\begin{proof}
For a vector $v = ^t(v_1, \dots ,v_{n-1}) \in P^{n-1}_n$, we set $\Supp (v) =
\{ i \mid v_i \neq 0\}$.
\begin{enumerate}
\item
If $\{ i,j \} \cap \Supp (v) = \emptyset$, then $\Supp (A_{ij}^{(n-1)}v) = 
\Supp (v)$ and $k$-th component of $A_{ij}^{(n-1)}v$ is equal to 
$a_{ij}v_k$ for $k \in \Supp (v)$.
\item
If $\{ i,j\} \cap \Supp (v) =\{ i \}$, then $\Supp (A_{ij}^{(n-1)}v)
\subset \Supp (v) \cup \{ j \}$ and $j$-th component and $i$-th component
of $A_{ij}^{(n-1)}v$ is equal to $-a_{nj}v_i$ and $(a_{ij} +
a_{nj})v_i$ respectively.  
\end{enumerate}
Therefore if we define $S_i$ inductively
by
$$
S_{i+1} = \begin{cases}
S_i & (\text{ if $e_{i+1} \cap S_i = \emptyset$}) \\
S_i \cup \{ j\} & (\text{ if $e_{i+1} \cap S_i = \{ k\}$ and $e_{i+1} = \{ j,k 
\}$}),  \\
\end{cases}
$$
and $S_0 = \{ p\}$.  
Then $\Supp (\prod_{i= l}^1A_{p_i, q_i}^{(n-1)})v
\subset S_{l}$.  If $q \in S_{l}$, then $t_1 < \cdots < t_m$.
This proves (1).  For (2), we can prove
$$
[(\prod_{i=s}^1A_{p_i, q_i}^{(n-1)}\left(\begin{matrix}
0 \\ \vdots \\ a_{nk} \\ \vdots \\ 0  \end{matrix}\right)
)]_{k_j} =
a_{n,k_j}\cdot
\prod_{j=s}^1 B_j \cdot
$$
if $t_j \leq s < t_{j+1}$  
by induction on $s$ using the infinitesimal pure braind relation (1) and (2).  
(In case $s \geq t_m$ 
and $s < t_1$, $a_{n,k_j} =a_{n,k_m}$ and $a_{n,k_j} = a_{n,k_0}$ respectively.)
This complete 2.
\end{proof}

Next we introduce an expression of $\emptyset ([r]) \wedge 
(r+1, i_{r+1}) \wedge \cdots \wedge (n, i_n)$ by using the notion of
principal graph.
\begin{definition}
For an index set $I = (i_{k+1}, \dots ,i_n)$, $(1 \leq i_p \leq p-1)$, we define
the ordered rooted graph $P_I$ as follows. 
\begin{enumerate}
\item
The set of vertices is $\{1, \dots ,n\}$, 
\item 
the set of roots is $\{ 1, \dots , k\}$, and
\item
the set of ordered edges is $\{(k+1, i_{k+1})< \dots < (n, i_n)\}$.  
\end{enumerate}
The graph $P_I$ is called the principal graph of $I$.  
\end{definition}
Let $p,q$ be two vertices contained in the
same connected component of $P_I$.  The unique shortest path connecting
$p, q$ in $P_I$ is denoted by $\gamma (p,q)$ 
and the minimal edge of $\gamma (p,q)$
is denoted by $\min (p,q)$.  Then by the construction of the 
principal graph, we have
the following lemma.
\begin{lemma}
Let us write a path $\gamma (p,q)$ connecting $p,q$ in $P_I$ as in figure 1:
Suppose that $e_{t_s}$ is the minimal edge of $\gamma (p,q)$.
Then $t_1 > \cdots > t_s < \cdots < t_m$
\end{lemma}

A graph $\Gamma$ is called a support of $\gamma = \sum_{\Gamma}a_{\Gamma}\Gamma$,
if $a_{\Gamma} \neq 0$.  The set of supports of $\gamma$ is denoted by
$\Supp (\gamma )$. 
Let $p,q \in [n]$ be vertices contained in the same connected
component in $P_I$.  We set $\gamma =\emptyset (\{1, \dots ,k\})\wedge
(k+1,i_{k+1}) \wedge \cdots \wedge (n, i_n)$.
By the construction of $\gamma$, if $\Gamma \in \Supp (\gamma )$ and $(p,q)$ 
appeares in $\Gamma$, then $(p,q)$ is the
$m$-th edge of $\Gamma$, where $e_m = \min (p,q)$.
Conversely, for any pairs $(p_{k+1},q_{k+1}),
\dots , (p_n, q_n)$ such that 
\begin{enumerate}
\item
$p_i$ and $q_i$ are contained in the same connected component of $P_I$, and
\item
$\min (p_j, q_j)$ is the $j$-th edge $(j,i_j)$ of $P_I$,
\end{enumerate}
$a_{\Gamma} = 1$ for $\Gamma = 
(p_{k+1},q_{k+1}) \cdots (p_n, q_n)$.  
We use distributive notation as
\begin{align*}
& \{(p_{k+1},q_{k+1})+(p'_{k+1},q'_{k+1})\}  
(p_{k+2},q_{k+2}) \cdots (p_n, q_n) \\
= &
(p_{k+1},q_{k+1})  (p_{k+2},q_{k+2}) \cdots (p_n, q_n)
+ (p'_{k+1},q'_{k+1})  (p_{k+2},q_{k+2}) \cdots (p_n, q_n). \\
\end{align*}
Here the right hand side has a meaning as a formanl linear combination of
ordered graphs.  The following proposition is nothing but the restatement of the
definition of $\wedge$.
\begin{proposition}
Let $S_i =\sum_{ 1 \leq p < q \leq n, \min (p,q) = e_i 
\text{ in } P_I} (p,q)$.  Then 
$$
\emptyset (\{1, \dots ,r\})\wedge (r+1, i_{r+1}) \wedge \cdots \wedge
(n,i_n) =
S_{r+1}\cdot S_{r+2} \cdots S_{n}
$$
\end{proposition}
We finish this subsection by computing 
$\Res_{x_n \to x_k}(\omega_{\Gamma})$ for an ordered rooted graph
$\Gamma = \{e_{r+1}, \dots ,e_{n}\} \in \Supp (\gamma )$.
Until the end of this subsection we assume
$\Gamma \in \Supp (\gamma )$ and $(n,k)$ is an edge of $\Gamma$.
Put $R_- = \{ k' \mid (n,k') \in \Gamma ,
\min (n,k') < \min (n,k)\}$ and $R_+ = \{ k' \mid (n,k') \in \Gamma ,
\min (n,k') \geq \min (n,k)\}$.
We meke a numbering of $R_+ =\{k_1 = i_n, k_2, \dots ,k_s = k\}$ such
that $\min(n,k_1) > \cdots > \min (n,k_s)$.  Set $e_{t_i}= \min (n,k_i)$.
For the figure of principal graph see figure 3.
If $s \geq 2$, we put $P = P(\Gamma , k)$ the power set of $R_+ -\{ k_1, k_s\}$.
For an element $p \in P$, we define a graph $\Gamma (p) \in \Gamma ([n-1],[r])$ 
as follows.
For $i = 2, \dots , s$, put $m(p,i) = \min 
\{j \mid k_j \in p\cup \{ k_1\}, j<i\}$.
The $t_i$-th edge of $\Gamma (p)$ is equal to $(k_i, k_m)$, where $m= m(p,i)$.
The $j$-th edge is the same as $\Gamma$ if $j \neq t_i, n$ $(i= 2, \dots , s)$.
A Set of ordered graph $\{ \Gamma (p) \mid p \in P(\Gamma ,k)\}$ is denoted 
by $R(\Gamma , k)$ and called the residue graph of $\Gamma$ with respect to
$k$.
\begin{figure}
%\label{fig:path}
\vskip 10pt
\centerline{\epsfbox{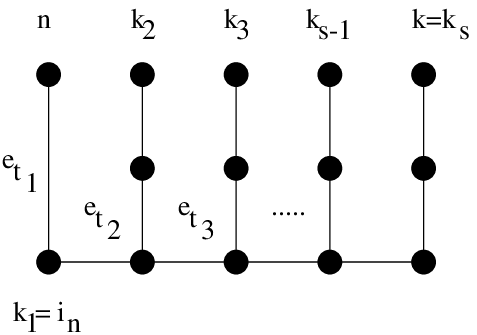}}
\vskip 5pt
\caption[]{ } 
\end{figure}
\begin{proposition}
\label{proposition:residue}
If $\# \mid R_+ \mid \geq 2$, then
$$
\Res_{x_n \to x_k}(\omega_{\Gamma}) =
\sum_{p \in P(\Gamma ,k)}(-1)^{\# p +1}\omega_{\Gamma (p)}.
$$
Here the residue $\Res_{x_n \to x_k}\omega = \eta \mid_{x_n = x_k}$,
where $\omega = d \log (x_n - x_k) \wedge \eta$.
\end{proposition}
\begin{proof}
We write $d \log (x_p - x_q) = \langle p,q \rangle$ for short.  
First we prove that
\begin{align*}
& \<< k, k_1 \>> \cdots \<< k, k_{s-1} \>>  \\
=&\sum_{p \subset \{ k_2, \dots , k_{s-1}\}} (-1)^{s+\# p}
\<< m(p,2),k_2 \>> \cdots \<< m(p,s-1),k_{s-1}\>>
\<< m(p,s), k_s\>> \\
\end{align*}
by induction on the cardinarity of $\{k_1, \dots ,k_{s-1}\}$.  
By the assumption of induction for 
$\{k_1, \dots k_{s-2}\}$, we have
\begin{align*}
&\<< k, k_1 \>> \cdots \<< k, k_{s-2} \>> \<<k, k_{s-1}\>> \\
= & \sum_{q \subset \{ k_2, \dots , k_{s-1}\}} (-1)^{s+\# q}
\<< m(q,2),k_2 \>> \cdots \<< m(q,s-1),k \>> \<< k,k_{s-1}\>> \\
= & \sum_{q \subset \{ k_2, \dots , k_{s-1}\}} (-1)^{s+\# q}
\<< m(q,2),k_2 \>> \cdots \<< m(q,s-2),k_{s-2} \>> \\
& \qquad (\<< m(q,s-1),k_{s-1} \>> \<< k,k_{s-1}\>> 
-\<< k_{s-1},m(q,s-1)\>> \<< k, m(q,s-1) \>>  \\
\end{align*}
and the last expression gives the expression of $\{k_1, \dots ,k_{s-1}\}$.
Therefore we have
\begin{align*}
& \Res_{x_n \to x_k} \<< n, k_1 \>> \cdots \<< n, k_s \>> \\
= &(-1)^{s-1}\<< n, k_1 \>> \cdots \<< n, k_{s-1} \>> \\
= & \sum_{p \subset \{ k_2, \dots , k_{s-1}\}} (-1)^{\# p+1}
\<< m(p,2),k_2 \>> \cdots \<< m(p,s-1),k_{s-1}\>>
\<< m(p,s), k_s\>>. \\
\end{align*}
This implies the proposition.
\end{proof}

\subsection{Proof of Theorem \ref{theorem:combinatorics}}
\label{subsec:proof of theorem}

We prove Theorem \ref{theorem:combinatorics} 
by induction. By Remark \ref{remark:homomorphism}, the morphism 
$$
\rho :A_{ij}^{(n-1)} \mapsto 
\Ind( \bold a )_{ij}
$$
defines a ring homomorphism from $P_{n-1}$ to $M(n-1, P_n)$.  We assume
Theorem \ref{theorem:combinatorics} for $n-1$.  
We define $W_{k} \in P_{n-1}^{r(r+1)\cdots (n-1)}$ for $k= 1, \dots ,n-3$ 
inductively by the relation similar to (\ref{eq:tensor vector}) and
$$
W_{n-2} = 
\left(\begin{matrix}
\frac{A^{(n-1)}_{n-1,1}}{x_{n-1}- x_1} \\
\vdots \\
\frac{A^{(n-1)}_{n-1,n-2}}{x_{n-1}- x_{n-2}} \\
\end{matrix}\right).
$$
Then by the assumption of induction, 
$$
\eta_{\gamma} = W_{r}(i_{r+1}, \dots ,i_{n-1})
dx_{n-1}\wedge \cdots \wedge dx_{r+1}
$$
for $\gamma = \emptyset (\{1, \dots ,r\}) \wedge (r+1, i_{r+1}) \wedge
\cdots \wedge (n-1, i_{n-1})$ in $P_{n-1}\otimes \Omega_R$.  
Here $W_r(i_{r+1}, \dots i_{n-1})$ is the $(i_{r+1}, \dots ,i_n)$-th
component of $W_r$.
By applying
the above ring homomorphism $\rho$, we have
$$
\rho(\eta_{\gamma}) =\rho ( W_{r}(i_{r+1}, \dots ,i_{n-1}))
dx_{n-1}\wedge \cdots \wedge dx_{r+1}
$$
in $M(n-1, P_{n-1}) \otimes \Omega_R$.  By the definition of $w_r$,  
$w_r (i_{r+1}, \dots , i_n)$ is equal to the $i_n$-th component of
the vector
$$
\rho ( W_{r}(i_{r+1}, \dots ,i_{n-1}))\left(\begin{matrix}
\frac{a_{n1}}{x_n -x_1} \\
\vdots \\
\frac{a_{nn-1}}{x_n- x_{n-1}} \\
\end{matrix}
\right).
$$
Therefore by taking the residue, $\Res_{x_n \to x_k}$, it is enough to
prove that
$$
\rho ( W_{r}(i_{r+1}, \dots ,i_{n-1}))\left(\begin{matrix}
0\\ \vdots \\
a_{nk} \\
\vdots \\
0 \\
\end{matrix}
\right)_{i_n} = \Res_{x_n \to x_k} (\eta_{\bar\gamma})
$$
for all $k = 1, \dots , n-1$,
where $\bar\gamma = \emptyset(\{1, \dots , k\})
\wedge (k+1, i_{k+1}) \wedge \cdots \wedge (n, i_n)$.  We compute
the left hand side and right hand side by using Lemma \ref{lemma:product}
and Proposition \ref{proposition:residue}.
Left hand side of (1) is expressed as a linear combination of $\eta_{\Gamma}$,
where $\Gamma$ is a support of $\gamma $. On the other hand, the expression
given in Proposition \ref{proposition:residue} gives an expression of the
right hand side by a linear combination
of $\eta_{\Gamma}$, where $\Gamma$ is a support of $\gamma $.
By comparing the coefficient of $\omega_{\Gamma}$
it is enough to prove the following 
proposition.
\begin{proposition}
Let $\Gamma \subset \Supp (\gamma)$ and $k, i_n$ be contained in the
same connected component of $\Gamma$. 
\begin{enumerate}
\item
$\length (k, i_n) = \# p +1$ if $\bar \Gamma (p) = \Gamma$.
Here $\length (k,i_n)$ is the length of the path connecting $k$ and
$i_n$ in $\Gamma$.
\item
$$
(-1)^{\length (k, i_n)}
\prod_{i= r+1}^n B_i = \sum_{\{\bar\Gamma \in \Supp (\bar\gamma )
\mid \Gamma \in R(\bar\Gamma, k) \}}
A_{\bar\Gamma}^{(n)},
$$
were $B_i$ is defined in Lemma \ref{lemma:product} and $R(\bar\Gamma , k)$
is defined in Proposition \ref{proposition:residue}.
\end{enumerate}
\end{proposition}
\begin{proof}
Let $\Gamma \in \Supp \gamma$ and suppose $R(\bar \Gamma,k) \ni \Gamma$.
As in Lemma \ref{lemma:product}, we make the numbering of the path 
in $\Gamma$ from $k$ to 
$i_n$ as $e_{t_1} = (n,k_1),
e_{t_2} = (k_1, k_2), \dots, e_{t_m} = (k_{s-1}, k)$ and $k_s = k, k_1 = i_n$.
First we claim the set $L= L(\bar\Gamma)=
\{ l \mid (n,l) \in \bar\Gamma \}$ contains 
$k_1, \dots ,k_s$.  Since $\Res_{x_n \to x_k }(\bar\Gamma )$
is not zero, $\bar\Gamma$ contains $(n,k)$, i.e. $k= k_s \in L$.
If the path connecting $k$ and $i_n$
in the corresponding graph $\bar\Gamma (p)$ is $k_1, \dots , k_s$,
then $p = \{ k_2, \dots ,k_{s-1}\}$.  Therefore $L \supset \{k_1, \dots ,k_s\}$.
If $l \in L -\{k_1, \dots , k_s\}$ and $q$ is a minimal element satisfying
$\min (i_n,l) < \min (i_n, k_q)$, then $\bar\Gamma (p)$ contins an edge 
$(l, k_q)$ by the definition of $\bar\Gamma (p)$, i.e. $(l, k_q) \in 
G(\Gamma ,k,i_n)$, where  
\begin{align*}
 G(\Gamma ,k,i_n)
 = \{e :\text{ edge } \mid & \text{ There exists $i$ such that }
e \ni k_i, \\
& \min (k_{i-1},i_n  ) \leq  e < \min (k_i, i_n)\}\}.
\end{align*}
Therefore $L-\{k_1, \dots ,k_s\} \subset G(\Gamma , k , i_n)$.

Conversely, for any subset $L$ of $\{ 1, \dots ,n\}$ satisfying 
\begin{enumerate}
\item
$L$ is contained in the same connected component of $i_n$,
\item
$L \supset \{k_1 ,\dots ,k_s\}$, and 
\item
$L - \{k_1, \dots ,k_s\} \subset G(\Gamma,k,i_n)$,
\end{enumerate}
there exists unique $\bar\Gamma (L)$ satisfying 
\begin{enumerate}
\item
$L (\bar \Gamma(L)) = L$, 
\item
$\bar\Gamma \in \Supp (\bar \gamma )$, and
\item
$\Supp (\Res_{x_n \to x_k}(\bar \Gamma)) \ni \Gamma$.
\end{enumerate}
Therefore
\begin{align*}
 \sum_{\{\bar\Gamma \mid \Gamma \in R(\bar\Gamma, k),
\bar\Gamma \in \Supp (\bar\gamma )\}}
A_{\bar\Gamma}^{(n)}  
= & \sum_{\substack{
L \supset \{k_1, \dots ,k_s\},
L-\{k_1, \dots, k_s\} \subset G(\Gamma,k,i_n),
\\
L \text{ is contained in the same} \\
\text{ connected component of }i_n}}
A_{\bar\Gamma (L)}^{(n)} \\
=(-1)^{\length (k, i_n)} & \prod_{i=k+1}^n B_i
\end{align*}
\end{proof}

\section{Proof of The Main Theorem}
\label{sec:proof of the main theorem}
\subsection{Some lemmata for the assypmtotic beheaviors}
\label{subsec:some lemmata}

In this subsection, we investigate the assymptotic beheavior of the
solution of linear differential equation with regular singularity.
Let $A \in \frac{1}{x}M(d, \Cal O_x)$, where $\Cal O_x$ is
a germ of holomorphic functions at $x=0$.  We are interested in 
the differential equation for $r\times r$-matrix valued function $V$:
$$
\frac{dV}{dx} = A V.
$$
We write $A = Rx^{-1} +\sum_{i=0}^{\infty}A_ix^i$, where $R, A_i \in
M(r,\bold C)$.  If all the eigen values of $R$ are enough small, then the
solution $V$ can be written as $V = F x^R C_0$, where $F$ is
an $r \times r$-valued homolorphic function in $I +x M(r,\Cal O_x)$, and 
$C_0 \in GL(r, \bold C)$.  In the rest of this section, we assume that
all the eigen values of $R$ are sufficiently small positive real numbers
and $R$ is semi-simple.  The eigen value of $R$ is denoted by
$0 < \lambda_1 <\cdots < \lambda_s$.
\begin{lemma}
\label{lemma:limit}
Let $\bold C^r = \oplus_{i=1}^s W_i$ be the eigen space decomposition of 
$\bold C^r$ with respect to $R$.
\begin{enumerate}
\item
If $w_i \in W_i$, then all the element $a_k$ of the vector
$Fx^Rw_i$ satisfies the estimation 
$\mid a_k \mid \leq \mid x \mid^{\lambda_i}c$ with some constant $c$ for
$k=1, \dots ,r$.  Moreover we have 
$\lim_{x\to 0}(x^{-\lambda_i}Fx^Rw_i) = w_i$.
\item
Let $\lambda > \lambda_i$. If $w_i \in W_i$ and all the elements $a_k$
of $Fx^Rw_i$ satisfy $\mid a_i \mid \leq x^\lambda c$ with some constant $c$,
then $w = 0$.
\item
Let $\lambda_i$ be an eigenvalue of $R$.  Let $p : W_i \to \bold C^l$
be a linear map and the composite $\bold C^r \to W_i \to \bold C^l$ is 
denoted by $\tilde p$.  Then we have
$$
\tilde p(\lim_{x \to 0}x^{-\lambda_i}F x^R w) =
\tilde p(\lim_{x \to 0}x^{-R} F x^{R}w)
$$
for any $w \in W$.
\end{enumerate}
\end{lemma}
\begin{proof}
Since $F = I + xm, m \in M(r, \Cal O_x)$, using identity 
$\lim_{x \to 0} x^{-\lambda}xm x^R = \lim_{x \to 0} x^{-R}xm x^{R} =0$,
we get the statements. 
\end{proof}
Let $n,k$ be integers such that $2 \leq k \leq n$ and we define
reduced part $V^{red} = V^{red}_k$ as in \S \ref{subsec:differential equation}.
The restriction of
$A_{ij}^{(k)}$ to $V^{red}$ is denoted by $A^{(k)}_{ij,red}$.  For
a subset $S$ of $[i,k]$, we define $A_S^{(k)}$ and $A_{S, red}^{(k)}$
by
$$
A_S^{(k)} =  \sum_{i<j, i,j \in S}A_{ij}^{(k)},
A_{S,red}^{(k)} =  \sum_{i<j, i,j \in S}A_{ij,red}^{(k)}. 
$$
Fron now on, $a_{ij}$ is sufficiently generic small positive real number.
For a semisimple matrix $A$, the formal sum of eigne values of $A$ counting
their multiplicities is denoted by $\sigma (A)$:
$\sigma (A) = \sum (\text{ eigen values of $A$})$.
In this situation, the set of eigen values is denoted by $\Supp (\sigma (A))$.
\begin{proposition}
\label{proposition:eigenvalue}
Under the notations and assumptions as above, $A_S^{(k)}$ and 
$A_{S,red}^{(k)}$ are semi-simple and 
\begin{align*}
\sigma (A_S^{(k)}) = & \sum_{T \subset [k+1,n]}
(k-l;\mid T^c \mid)(l ; \mid T \mid) a_{S\cup T}, \\
\sigma (A_{S,red}^{(k)}) = & \sum_{T \subset [k+1,n]}
(k-l-1;\mid T^c \mid)(l ; \mid T \mid) a_{S\cup T}, \\
\end{align*}
where $a_U = \sum_{i< j, i,j \in U} a_{ij}$ for a subset $U \subset [1,n]$.
For a subset $T \subset [k+1,n]$, $T^c = [k+1, n] -T$ and 
$l = \#\mid S \mid -1$ and $(a;b) = a(a+1) \cdots (a+b-1)$.
\end{proposition}
To prove the above propostion, we use the following two elementary
lemmata.
\begin{lemma}
Let $X$ be $kN \times kN$-matrix. We assume that there exist
semi-simple matrices $B$ and $D$ and matrices $C_1, \dots ,C_k$
such that
\begin{align*}
A\left(\begin{matrix}
0  \\
\vdots  \\
1 \\
-1  \\
\vdots   \\
0  \\
\end{matrix}
\right) 
\begin{matrix}
 \\
 \\
 i \\
 i+1 \\
 \\
 \\
\end{matrix}
& =
\left(\begin{matrix}
0  \\
\vdots  \\
B \\
-B  \\
\vdots   \\
0  \\
\end{matrix}
\right)  \\
A\left(\begin{matrix}
C_1  \\
\vdots  \\
C_k \\
\end{matrix}
\right) & =
\left(\begin{matrix}
C_1D  \\
\vdots  \\
C_kD \\
\end{matrix}
\right),
\end{align*} \\
with $\Supp (\sigma (B)) \cap \Supp (\sigma (D)) = \emptyset$.
Then 
\begin{enumerate}
\item
$\sigma (A) = (k-1)\sigma (B) + \sigma (D)$.
\item
$(k-1)N$-dimensional subvector space $V^{red} =
\{ (v_1, \dots ,v_k) \mid v_i \in \bold C^N, \sum v_i = 0\}$
is stalbe under the action of $A$.  Let $A^{red}$ be the restriction
of $A$ to $V^{red}$. Then $\sigma (A^{red}) = (k-1)\sigma (B)$.
\end{enumerate}
\end{lemma}
\begin{lemma}
Let $a_{ij} \in P_k$ and set 
$A_{ij}= \Ind (\bold a )_{ij}$
for $1 \leq i < j \leq k-1$,  
$A_{[1,k-1]} = \sum_{1 \leq i < j \leq k-1}A_{ij}$,
$a_{[1,k-1]} = \sum_{1 \leq i < j \leq k-1}a_{ij}$ and
$a_{[1,k]} = \sum_{1 \leq i < j \leq k}a_{ij}$.
Then we have
\begin{align*}
A_{[1, k-1]}\left(\begin{matrix}
0  \\
\vdots  \\
1 \\
-1   \\
\vdots   \\
0  \\
\end{matrix}
\right) 
\begin{matrix}
 \\
  \\
i \\
i+1   \\
   \\
  \\
\end{matrix}
& =
\left(\begin{matrix}
0  \\
\vdots  \\
a_{[1,k]} \\
-a_{[1,k]} \\
\vdots   \\
0  \\
\end{matrix}
\right)  \\
A_{[1,k-1]}\left(\begin{matrix}
a_{k1}  \\
\vdots  \\
a_{k-1} \\
\end{matrix}
\right) & =
\left(\begin{matrix}
a_{k1}a_{[1,k-1]}  \\
\vdots  \\
a_{kk-1}a_{[1,k-1]}  \\
\end{matrix}
\right).
\end{align*} \\
\end{lemma}
\begin{proof}
The first equality follows from the expression
$$
A_{[1,k-1]}=
\left(\begin{matrix}
a_{[1,k-1]}+\sum_{j \neq 1} a_{k,j} & -a_{k1} & \cdots \\
-a_{k2} & a_{[1,k-1]} + \sum_{j \neq 2}a_{kj} & \cdots \\
-a_{k3} & -a_{k3} & \cdots \\
\vdots & & & \vdots \\
\end{matrix}\right).
$$
The second equality is obtined directly by the equality
$$
A_{ij}\left(\begin{matrix}
a_{k1} \\ \vdots \\ a_{kk-1} \\
\end{matrix}\right) =
\left(\begin{matrix}
a_{k1}a_{ij} \\ \vdots \\ a_{kk-1}a_{ij}
\end{matrix}\right).
$$
\end{proof}
\begin{proof}(of Proposition \ref{proposition:eigenvalue})  
We prove the proposition by induction.
By two lemmata, we have
\begin{align*}
\sigma (A_S^{(k)}) & = (k-l) \sigma (A_S^{(k+1)}) + l
\sigma (A_{S\cup\{k+1\}}^{(k+1)}), \\
\sigma (A_{S,red}^{(k)}) & = (k-l-1) \sigma (A_{S,red}^{(k+1)}) + l
\sigma (A_{S\cup\{k+1\},red}^{(k+1)}), \\
\end{align*}
using homomorphism $P^{(k+1)} \to M((k+1)(k+2) \cdots (n-1), \bold C)$
and assumption of independence of $a_{ij}$.
\end{proof}

\subsection{Relation between Selberg integral and Drinfeld associator}
\label{subsec:relation between selberg integral}
In this section, we will compare vectors whose elements are given by 
Selberg integrals with the Drinfeld associator.  Let $n \geq 3$ be an
integer and we define $A_{ij}^{(k)}$ as in \S \ref{subsec:differential equation}.
We set 
$V = V_{3,n} = \bold C^{3\cdot 4 \cdots (n-1)}$.  
Let 
$S = S([n]/[3],x_1, x_2, x_3, \alpha_{ij})$ be a $V$-valued function on $x_1, x_2, x_3$
whose $(i_4, \dots , i_n)$-component is given by 
$S_{\emptyset(\{1,2,3\})\wedge (4,i_4) \cdots (n,i_n)}
([n]/[3],x_1, x_2, x_3, \alpha_{ij})$.  
Then $S$ satisfies the differential equation
$$
dS = (A_{13}^{(3)}d \log (x_1 - x_3) +
A_{23}^{(3)}d \log (x_2 - x_3) ) S.
$$
We set $\bar S(x_3) = S([n]/[3],0,1,x_3)$.  Then $\bar S$ satisfies the equation
$$
\frac{d\bar S}{dx_3} =
(A_{13}^{(3)}\frac{dx_3}{x_3} +
A_{23}^{(3)}\frac{dx_3}{x_3-1} )\bar S.
$$
Therefore by considering the rational representation $\rho$ of degree 1 :
$\rho : \bold C \langle\langle  X, Y \rangle\rangle  \to 
M(3\cdot 4 \cdots (n-1), \bold Q[[\alpha_{ij}]])$, we have
$$
\lim_{x \to 1}(1-x_3)^{-A_{23}^{(3)}} \bar S(x_3) = 
\rho (\Phi (X,Y)) \lim_{x_3 \to 0}
x_3^{-A_{13}^{(3)}}\bar S(x_3).
$$
We have the following lemma.
\begin{lemma}
\label{lemma:limit to 0}
\begin{enumerate}
\item
For any $i_4, \dots ,i_n$, we put $\gamma = \emptyset(\{1,2,3\})\wedge
(4, i_4)\wedge \cdots \wedge(n,i_n)$.  Then for a sufficiently small $x_3$,
we have an estimation 
\begin{align}
\label{eq:estimate at 0}
\mid S_{\gamma}([n]/[3],0,1,x_3)\mid  < c x_3^{\alpha_{max}}
\end{align}
for some constant $c$.  Here $\alpha_{max}$ is the maximal eigenvalue
\linebreak
$\sum_{1 \leq i < j \leq n, i,j \neq 2}\alpha_{ij}$
of $A_{13}^{(3)}$.
\item
For $\Gamma \in \Gamma ([n],[3])$,
\begin{align*}
& \lim_{x_3 \to 0}x_3^{-\alpha_{max}}
S_{\Gamma}([n]/[3], 0,1,x_3)
\\
= &
\begin{cases}
S_{\Gamma '}([n]-\{ 2\}/\{ 1, 3\}, 0,1)
 &
(\text{if  there is no edges containing 2}) \\
0 & (otherwise) \\
\end{cases} \\
\end{align*}
Here $\Gamma '\in \Gamma ([n]-\{ 2 \},\{1, 3\})$ 
is the ordered graph obtained by deleting $2$ from
the graph $\Gamma$.
\end{enumerate}
\end{lemma}
\begin{proof}
By Proposition \ref{proposition:eigenvalue}, we have $\alpha_{max} =
\sum_{1 \leq i < j \leq n, i,j \neq 2} \alpha_{ij}$.  To prove
the statement, it is enough to prove that
$$
\int_D \prod_{1 \leq i < j \leq n}(x_i -x_j)^{\alpha_{ij}}
\omega_{\Gamma} \mid_{x_1 = 0, x_2=1}
$$
satisfies the estmation of (\ref{eq:estimate at 0}) for an ordered
rooted tree $\Gamma$ with the root set $[3]$.  We change variable
by $x_p = \xi_p x_3$ for $p = 4, \dots ,n$.  Then
\begin{align}
\label{eq:coordinate change}
\omega_{\gamma} = \pm \prod_{(p_i, q_i) \in E_{\Gamma},\text{
not adjacent to 2}}\frac{d\xi_{p_i}-d\xi_{q_i}}{\xi_{p_i}-\xi_{q_i}}
\prod_{(p_i, 2) \in E_{\Gamma}}\frac{x_3 d\xi_{p_i}}{-1}
\cdot (1+o(1)).
\end{align}
and 
$$
\prod_{1 \leq i < j \leq n}(x_i - x_j)^{\alpha_{ij}}
= 
\prod_{1 \leq i < j \leq n, i,j \neq 2}(\xi_i -\xi_j)^{\alpha_{ij}}
\cdot x_3^{\alpha_{max}} \dot (1+ o(1)).
$$
Here we put $\xi_3 =1, \xi_1 =0$.
In particular, $\lim_{x_3 \to 0}x_3^{-\alpha_{max}}\int_D
\Phi \omega_{\Gamma} = 0$ if $\Gamma$ contains an edge adjacent to 2.
The signature in (\ref{eq:coordinate change}) arise from
the subtutition for separating edges of $\Gamma$ adjacent to 2 and those does not
adjacent to 2.
If $\Gamma$ contains no edges adjacent to 2, we get the second statement.
\end{proof}
From Lemma \ref{lemma:limit}, we have the following corollary.
\begin{corollary}
The $(i_4, \dots ,i_n)$-th component of 
$\lim_{x_3 \to 0}x_3^{-A_{13}^{(3)}}\bar S(x_3)$ is equal to 
$S_{\gamma}([n]-\{ 2\}/\{1, 3\},0,1)$
if $i_p \neq 2$ for $p= 4, \dots ,n$, where 
$\gamma = \emptyset(\{1,3\})\wedge (4, i_4) \wedge \cdots \wedge (n, i_n)$
and $0$ otherwise.
\end{corollary}
\begin{proof}
By the definition of $\gamma = \emptyset ([3])\wedge (4, i_4) \wedge
\cdots \wedge (n, i_n)$, if $i_p=2$ for some $p$, then any
$\Gamma \in \Supp (\gamma )$ has an edge adjacent to 2.
If $i_p \neq 2$ for all $p$, then any $\Gamma \in \Supp (\gamma )$ contains
no edges adjacent to 2.  Therefore the statement follows from Proposition 
\ref{lemma:limit to 0}.
\end{proof}
Next we consider the assumptotic beheavior for $x_3 \to 1$.  Let
$I$ be the set $\{ (i_4, \dots ,i_n) \mid i_p \neq 2,3\}$.  By the definition
of $A_{23}^{(3)}$, the projection $p : V \to \bold C^I$ to $I$-th coordinate
factors through $\alpha_{23}$ eigen projection.
Therefore we have
$$
p( \lim_{x_3 \to 1}(1-x_3)^{-A_{23}^{(3)}}\bar S(x_3)) =
p(\lim_{x_3 \to 1}(1-x_3)^{-\alpha_{23}}\bar S (x_3)).
$$
by Lemma \ref{lemma:limit}.
On the other hand, it is easy to see the following lemma.
\begin{lemma}
\label{lemma:limit to 1}
If $\Gamma$ contains no edges containing 2 and 3, then
$$
\lim_{x_3 \to 1}(1-x_3)^{-\alpha_{23}} 
S_{\Gamma}([n]/[3], 0,1,x_3)
 =
S_{\Gamma '}([n]-\{ 3\}/[2], 0,1,\alpha_{ij}')
$$
where $\Gamma '$ is the ordered graph obtained by deleting 3 from the graph
$\Gamma$, and 
$\alpha_{ij} = \alpha_{ij}'$ if $i,j \neq 2$ and 
$\alpha_{2j}' = \alpha_{2j}+\alpha_{3j}$.
\end{lemma}
\begin{definition}
The vector $\lim_{x_3 \to 0} x_3^{-A_{13}^{(3)}}\bar S(x_3)$
and $\lim_{x_3 \to 1} (1-x_3)^{-A_{23}^{(3)}}\bar S(x_3)$ is denoted by
$V^{(1)}$ and $V^{(2)}$ respectively.
Then we have
\begin{align}
\label{eq:projection}
p(V^{(2)}) = p(\rho (\Phi (X,Y))V^{(1)})
\end{align}
\end{definition}
By Lemma \ref{lemma:limit to 1}, the $(i_4, \dots ,i_n)$-th component of
$V^{(2)}$ with $i_p \neq 2,3$
is equal to $S_{\gamma}(0,1, \alpha'_{ij})$,
where $\alpha'_{ij} = \alpha_{ij}$ if $i,j \neq 2$ and $\alpha '_{2,j} =
\alpha_{2j} +\alpha_{3j}$, where $\gamma = \emptyset (\{1,2\})
\wedge (4,i_4) \wedge \cdots \wedge (n, i_n)$.
We compute the limit of all the component of $V^{(1)}$ for the limit
$\alpha_{3i} \to 0$. For this purpose, we compute 
$\lim_{\alpha_{3i}\to 0}S_{\gamma}(\alpha_{ij})$ for 
$\gamma = \emptyset (\{1,3\})\wedge (4, i_4) \wedge \cdots \wedge (n,i_n)$
with $i_p \neq 2$, in the next subsection.

\subsection{Limit for $\alpha_{3i} \to 0$}
\label{sebsec: limit for a}

In this subsection, we change numbering from that of the last subsection.
Let $\Gamma$
be an ordered graph with the root set $[2]$ and vertex set $[n]$.
We set $\Phi = \prod_{i\ll j}(x_i - x_j)^{\alpha_{ij}}$, and
$$
S(\alpha_{ij})=
\int_{D([n]/[2],0,1)}\eta_{\Gamma}\Phi .
$$
Before proving Proposition \ref{proposition:alpha limit},
we remark the following lemma.
\begin{lemma}
Let $F(x)$ be a continuous function defined on $(p,1]$.
Suppose that $F(x)$ is integralble on $(p, p+\epsilon]$.  Then we have
$$
\lim_{\alpha \to 0}\int_p^1 \alpha (1-x)^{\alpha -1}F(x) dx = F(1).
$$
\end{lemma}
\begin{proof}
This is the fundamental property of $\delta$-function 
$\lim_{\alpha\to 0}\alpha (1-x)^{\alpha-1}$.
\end{proof}
\begin{proposition}
\label{proposition:alpha limit}
\begin{enumerate}
\item
If $\lim_{\alpha_{2i} \to 0}S(\alpha_{ij}) \neq 0$, then
(1) $\Gamma$ contins no edges adjacent to 2, or (2)
$(2, 3)$ is a unique edge adjacent to 2.
\item
If $(2, 3)$ is a unique edge in $\Gamma$ adjacent to 2,
then $\lim_{\alpha_{2i} \to 0}S(\alpha_{ij})$ is
equal to $S_{\Gamma '}(\alpha_{ij}')$, 
where $\Gamma '$ is the ordered
graph obtained by deleting the edge $(2,3)$ and replace the numbering
$3$ of original edge by the new numbering 2 and $\alpha_{ij}' = \alpha_{ij}$
if $i,j \neq 2$ and $\alpha_{2,k}' = \alpha_{3,k}$.
\end{enumerate}
\end{proposition}
\begin{proof}
Suppose that $\Gamma$ contains an edge adjacent to 2.  Let $p \leq 3$
be the minimal number such that $(2,p)$ is an edge of $\Gamma$.  
Set 
\begin{align*}
 F(x_p, \dots ,x_n) 
= & 
\prod_{(pq)\in E_{\Gamma}, \neq (2,p)}a_{pq}
\prod_{1 \leq i\leq n, p \leq j \leq n, i < j, (i,j) \neq (2,p)}
(x_i-x_j)^{\alpha_{ij}+\epsilon_{ij}} \\
& \int_{\{x_p < \cdots < x_3 < 1\}}
\prod_{1 \leq i < j \leq p -1}(x_i -x_j)^{\alpha_{ij}+ \epsilon_{ij}}
dx_{p-1} \cdots dx_3,\\
\end{align*}
where $\epsilon_{ij}=-1$ if $(i,j)$ is an edge of $\Gamma$ and $0$
otherwise.
Then
$$
\lim_{\alpha_{2p} \to 0}
S_{\Gamma} =
\int_{\{0 < x_n < \cdot < x_p <1\}}
F(x_p, \dots ,x_n)\alpha_{2p}(1-x_p)^{\alpha_{2p}-1}.
$$
Therefore if $p\neq 3$ or there exist at least tow p such that
$(2,p)$ is an edge of $\Gamma$, then $S_{\Gamma} = 0$. If $p=3$ and
there is no edge adjacent to 2 other than $(2,3)$, 
then
$\lim_{\alpha_{23} \to 0}S_{\Gamma} = S_{\Gamma '}(\alpha_{ij}')$.
\end{proof}
We define $S_{\gamma}(\alpha_{ij})$ by $\sum a_{\Gamma}S_{\Gamma}(\alpha_{ij})$,
where $\gamma = \sum a_{\Gamma}\Gamma
\in \Gamma ([2],[n])$.
\begin{corollary}
\label{corollary:alpha limit}
Let $\gamma = \emptyset (\{1,2\})\wedge (3, i_3) \wedge \cdots \wedge
(n, i_n)$.  
\begin{enumerate}
\item
If there exists $k \neq 3$ such that $i_k =2$, then
$\lim_{\alpha_{2i} \to 0}S_{\gamma}(\alpha_{ij}) =0$
\item
If $i_3 = 2$ and $i_k \neq 2$ for $k \neq 3$, then 
$$
\lim_{\alpha_{2i} \to 0}S_{\gamma}(\alpha_{ij})=
S_{\gamma '}(\alpha_{ij}'),
$$
where $\gamma '$ is $\emptyset (\{1,3 \}) \wedge (4, i_4) \wedge \cdots \wedge
(n, i_n)$.
\end{enumerate}
\end{corollary}
\begin{proof} (Proof of the Main Theorem \ref{theorem:the main theorem})
We can proceed by the induction on $n$.
We consider the limit of ($\ref{eq:projection}$) for $\alpha_{3i} \to 0$.  Then
all the entries of  
$\lim_{\alpha_{3i} \to 0}(\rho (\Phi (X,Y)))$ are contined in $H_{\alpha}$ by
Corollary \ref{corollary:homogeneous rational representation}.
By Corollary \ref{corollary:alpha limit}, all the entries of
$\lim_{\alpha_{3i} \to 0}V^{(1)}$ are contained in $H_{\alpha}$. Therefore
all the entries of $\lim_{\alpha_{3i} \to 0}p(V^{(0)})$ are also contained in
$H_{\alpha}$. Therefore $S_{\gamma}(0,1, \alpha_{ij})$ is an element of 
$H_{\alpha}$ for $\gamma = \emptyset (\{1, 2\}) \wedge (3, i_3) \wedge
\cdots \wedge (n, i_n)$ under the restriction 
\begin{quote}
(R) : $i_k \neq 2$ for all $k$.  
\end{quote}
On the other hand,
by the relation (\ref{eq:fundamental relation}), 
the restriction (R) is not necessary.
This completes the main therem.
\end{proof}

\end{document}